\makeatletter
\AddToHook{class/l4dc2027/after}{%
  \let\hl@conference@titlepagestyle\ps@jmlrtps
  \def\ps@jmlrtps{%
    \hl@conference@titlepagestyle
    \let\@oddhead\@empty
    \let\@evenhead\@empty
  }%
}
\makeatother
\PassOptionsToPackage{numbers,square,sort&compress}{natbib}
\documentclass{l4dc2027}
\usepackage{iftex}
\ifXeTeX
  \usepackage{fontspec}
\else
  \usepackage{times}
\fi

\usepackage{booktabs}
\usepackage{array}
\usepackage{graphicx}
\usepackage{tikz}
\usepackage{amsmath,amssymb}
\usepackage{xcolor}
\usepackage{microtype}
\bibpunct{[}{]}{,}{n}{,}{,}  

\title[ControlHarbor: A Benchmark for Control Design by Coding Agents]{ControlHarbor: A Benchmark for Control Design by Coding Agents}

\author{\Name{Jieqiang Wei} \Email{jieqiang.wei@eqtpartners.com}\\
        \addr EQT, Stockholm, Sweden
        \AND
        \Name{Fahad Ali} \Email{fahad@redpine.ai}\\
        \addr Redpine, Stockholm, Sweden
        \AND
        \Name{Jonny Reichwald} \Email{jonny.reichwald@eqtpartners.com}\\
        \addr EQT, Stockholm, Sweden
        \AND
        \Name{Silun Zhang} \Email{silunz@kth.se}\\
        \addr Department of Mathematics, KTH Royal Institute of Technology, Stockholm, Sweden}

\begin{document}

\maketitle

\begin{abstract}
LLM agents, i.e., language models run by a harness, can now write, run and debug software autonomously.
This raises the question of whether they can also carry out control system design. In this paper, we
consider control system design as a task for LLM agents and grade the designed controller as control
engineers do: the controller is given in a single Python file mapping measurements to actions,
evaluated on held-out seeds against an ordered set of requirements, under specified physical
perturbations and structural certificates, without using neural networks. We formulate the task as a
chance-constrained feasibility problem and instantiate it in ControlHarbor, a benchmark of 44 tasks over
13 plants, which includes network consensus, vehicle platooning, the MuJoCo hopper, walker and humanoid.
The tasks are packaged in the Harbor format, so that any agent harness can be run against them, and
the same environments also serve learning algorithms other than LLMs, such as reinforcement learning
(RL), under the same evaluation. With ControlHarbor, we observe that the controllers designed by three
frontier coding agents, Claude Opus 5.5, Fable 5.1 and Sonnet 5.5, reach the RL return on every parity
task, but only in a minority of sessions. We also observe that the controllers certified for return
alone usually fail under perturbations. When robustness is part of the specification, a structured
design meets it at the RL return on the humanoid task, where domain-randomised RL does not; under large
pushes, domain-randomised RL succeeds and no designed controller does. The benchmark, the controllers,
and their grades are released.\footnote{Tasks, controllers, grades and scripts:
\url{https://github.com/jq-wei/hl-control}.}
\end{abstract}


\section{Introduction}
\label{sec:intro}

Coding agents, i.e., language models operating a terminal, now solve a growing share of the
software-engineering tasks of benchmarks such as SWE-bench and Terminal-Bench
\citep{jimenez2024swebench,terminalbench2026}, and the harnesses built for those benchmarks, i.e., the
programs that call the model in a loop, give it tools and record its actions, run any such
agent in a container against a test it cannot see \citep{harbor2026}. Control design is also done in
code, as a feedback law and its evaluation in simulation, so a natural question is \textbf{whether LLM
agents can design controllers, and how a controller designed this way compares with a
reinforcement-learning solution, trained on the same plant}. In this paper we consider control design as a task for a
coding agent and grade it as control engineers do: one file, evaluated on seeds the agent never sees
against an ordered set of requirements, under physical perturbations and structural certificates,
within one hour, and without using neural networks or gradient-based learning.\footnote{Both are excluded
so that the controller stays a white-box program whose structure an engineer can read, check and
modify; the RL solutions supply the black-box comparison.} The controller is given in a single
Python file that maps measurements to actions: the file is the whole artefact, so it can be read,
re-run and graded as it stands, and it cannot rely on trained weights or outside services. For each
plant, successive task tiers introduce one additional challenge at a time, which are model mismatch,
sensing delay and noise, heterogeneity, coupling and disturbances. On the hopper, walker, humanoid, and
lunar lander, the designed controllers are compared with the reference RL policies under the same
plant-specific perturbations on the robust tiers. On the three locomotion plants, we also test nominal
parity by requiring controllers to meet a distance target and achieve at least 90\,\% of the reference
policy's mean return. Together, these tests show which performance and robustness requirements each
approach meets on each plant. The locomotion plants serve as stress tests: they are hard, far from saturated,
and carry strong RL solutions, and they are not a model of how legged robots are controlled in
deployment.

\begin{figure}[t]
\centering
\newcommand{\plantpanel}[2]{\begin{minipage}[t]{0.24\textwidth}\centering
  #1\\[-1pt]{\small #2}\end{minipage}}
\plantpanel{\includegraphics[width=\linewidth]{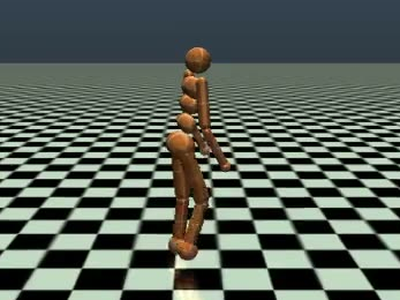}}{(a) Humanoid}\hfill
\plantpanel{\includegraphics[width=\linewidth]{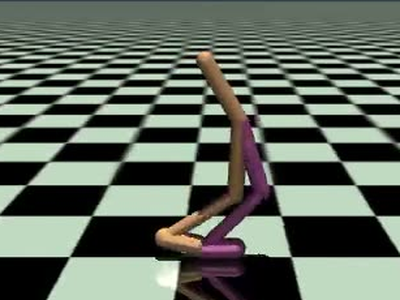}}{(b) Walker}\hfill
\plantpanel{\includegraphics[width=\linewidth]{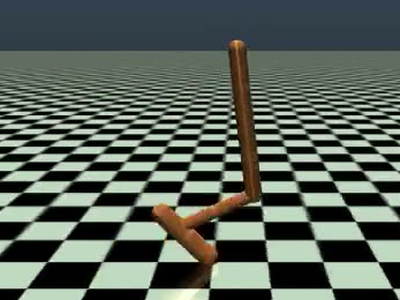}}{(c) Hopper}\hfill
\plantpanel{\includegraphics[width=\linewidth]{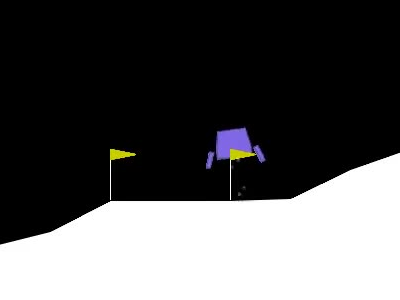}}{(d) Lunar lander}\\[4pt]
\plantpanel{\includegraphics[width=\linewidth]{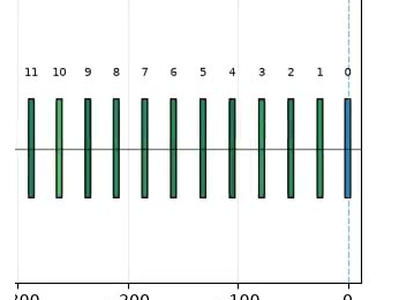}}{(e) Vehicle platoon}\hfill
\plantpanel{\includegraphics[width=\linewidth]{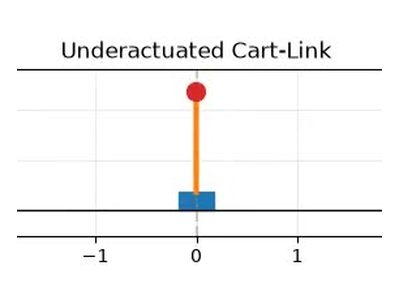}}{(f) Cart-link swing-up}\hfill
\plantpanel{\includegraphics[width=\linewidth]{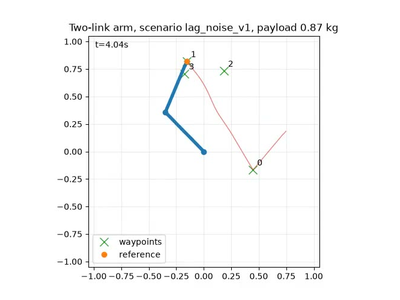}}{(g) Two-link arm}\hfill
\plantpanel{\includegraphics[width=\linewidth]{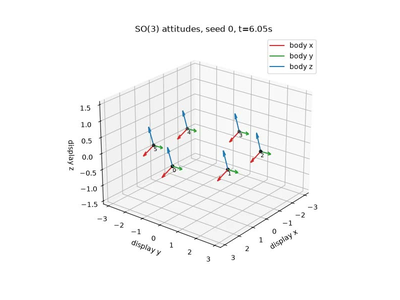}}{(h) SO(3) synchronisation}
\caption{The benchmark spans networked systems, manipulators and legged robots, and every plant shown is
driven by a controller an agent wrote in one session. Eight of the 13 plants:
(a)--(c) the MuJoCo humanoid, walker and hopper; (d) the Box2D lunar lander; (e) the lead vehicles of a
platoon of 20 with delayed, noisy local sensing; (f) the cart-link swing-up on a bounded track;
(g) a two-link arm with actuator lag and an unknown payload; (h) six rigid bodies synchronising their
attitude on SO(3) over a ring.}
\label{fig:plants}
\end{figure}

Related work falls into three groups, by what the model is asked to produce. In the first, it answers
or tunes textbook control: ControlBench \citep{kevian2024controlbench} grades frontier models on
undergraduate problems, ControlAgent and AgenticControl \citep{guo2024controlagent,narimani2025agenticcontrol}
tune classical compensators on linear plants, and CoDyControlBench \citep{zhou2026codycontrolbench} asks
six models to design feedback for 132 coupled systems; these check the specification on the nominal
model and put no RL solution beside the designed one. In the second, the model writes an executable
policy and a search loop scores it: Code as Policies and Voyager \citep{liang2022code,wang2023voyager}
write robot programs, Eureka \citep{ma2024eureka} writes the reward, and FunSearch, AlphaEvolve and GEPA
\citep{romera2024funsearch,novikov2025alphaevolve,agrawal2026gepa} search over programs or prompts with a
scorer in the loop. With a coding agent maintaining a programmatic policy, the pattern has been called
heuristic learning \citep{weng2026beyondgradients,wei2026blog}; it has been applied to the lunar lander
\citep{guo2026codeevolution}, to steel rolling \citep{siboni2026steel}, to manipulation
\citep{elmaaroufi2026rho} and, independently, to active flow control, where it matches or beats deep
reinforcement learning on 10 of 13 benchmarks \citep{garnier2026flow}. Closest to us,
\citet{ahmed2026codetocontrol} let the model write the controller's structure and a derivative-free
search fit its parameters, reaching PPO-level return on MuJoCo. Each of these works uses its own search
loop and its own scorer, compares with reinforcement learning on its own terms, tests robustness after the fact if at
all, and offers no task format that other agents can be run against. Closer to deployment, language-model
agents now act as the policy itself on generalist manipulation benchmarks such as RoboDojo
\citep{robodojo2026}. In one harness the agent calls motion, perception and state tools from Python, with
inverse kinematics and motion execution supplied below that interface, and the remaining failures are
dominated by grasps that slip, which originate below the interface \citep{hexashell2026}; the developer
of Claude Fable 5.2 (released on 2026-11-09) reports an average success of 64.9\,\% in simulation, against 27.6\,\% for Fable 5.1,
from its own internal runs \citep{anthropic2026fable52}. In these works the low-level feedback is
supplied by the platform; here it is the object the agent is asked to write. The third group grades RL solutions and designed controllers on the same plants: the MuJoCo suite and its public RL solutions
\citep{todorov2012mujoco,towers2024gymnasium,raffin2021sb3,haarnoja2018sac,fujimoto2018td3,schulman2017ppo},
the linear policies of \citet{mania2018ars}, to which our results return, the control view of
reinforcement learning \citep{recht2019tour} and of data-driven design \citep{dorfler2023datadriven},
against the designed side of SIMBICON, Raibert hopping and capture-point walking
\citep{yin2007simbicon,raibert1986legged,pratt2006capture}; the two sides are rarely graded by one gate
on one set of seeds.

We make three contributions. First, we formulate controller design as a chance-constrained feasibility
problem and use held-out seeds to test whether a controller meets the specified requirements on nominal
and perturbed plants (Section~\ref{sec:grader}). Second, we introduce ControlHarbor, a benchmark of 44
tasks across 13 plants (Figure~\ref{fig:plants}), and release the tasks, the passing controllers and
their grades, and the scripts used to generate the tables (Appendices~\ref{app:harbor}
and~\ref{app:repro})\footnote{The rendered attempts are at \url{https://jq-wei.github.io/hl-control/}.}.
Third, we evaluate three frontier coding models across 330 sessions (Section~\ref{sec:boundary},
Table~\ref{tab:boundary}). Agents met the nominal parity criteria---a distance target and at least
90\,\% of the reference RL policy's mean return---in 1/9 hopper, 3/9 walker, and 3/9 humanoid sessions.
All three humanoid parity controllers failed under weaker motors. When the task also required
robustness, 2/9 sessions produced SIMBICON controllers that passed the full suite, a target missed by
SAC trained with randomised physics and pushes. A second SAC policy trained with additional noise and
delay passed the separate 100~N push test that all nine agent sessions failed.

\begin{table}[t]
\centering\scriptsize
\setlength{\tabcolsep}{3pt}
\begin{tabular}{@{}>{\raggedright\arraybackslash}p{2.65cm}>{\raggedright\arraybackslash}p{2.85cm}>{\raggedright\arraybackslash}p{3.9cm}l>{\raggedright\arraybackslash}p{3.2cm}@{}}
\toprule
task & asked for & RL solution & O / F / S & best designed in session \\
\midrule
\addlinespace[2pt]\multicolumn{5}{@{}l}{\emph{Hopper}} \\
\texttt{forward} & forward 3 m & TD3 3532, 20.3 m: clears & 3/3 / 3/3 / 2/4 & return 2605 (O, Raibert) \\
\texttt{parity} & fast 18 m, return $\ge$ 3179 & TD3 3532 (the bound is TD3 $-$ 10\,\%) & 1/3 / 0/3 / 0/3 & return 3450 (O, linear policy (search)) \\
\texttt{robust} & forward under 4 perturbations & TD3 alive 0.62 / 0.00 / 0.00 / 1.00: fails & 3/3 / 2/3 / 1/3 & return 2145 (F, other) \\
\texttt{brisk-robust} & brisk 10 m under the suite, return $\ge$ 2500 & TD3 fails; SAC-DR fails & 1/3 / 0/3 / 2/3 & return 2652 (O, other) \\
\addlinespace[2pt]\multicolumn{5}{@{}l}{\emph{Walker2d}} \\
\texttt{forward} & forward 3 m & TD3 4701, 29.6 m: clears & 3/3 / 3/3 / 3/3 & return 2494 (S, CPG) \\
\texttt{parity} & fast 20 m, return $\ge$ 4231 & TD3 4701 (bound TD3 $-$ 10\,\%) & 2/3 / 1/3 / 0/3 & return 5265 (O, linear policy (search)) \\
\texttt{robust} & forward under 4 perturbations & TD3 falls at motors $-20\,\%$: fails; SAC-DR clears & 3/3 / 3/3 / 3/3 & return 2588 (F, SIMBICON) \\
\addlinespace[2pt]\multicolumn{5}{@{}l}{\emph{Humanoid}} \\
\texttt{forward} & forward 3 m & SAC 5579: clears & 2/3 / 2/3 / 1/3 & return 5522 (F, Raibert) \\
\texttt{brisk} & brisk 8 m & SAC brisk 0.975: clears & 2/3 / 1/3 / 1/3 & return 5761 (O, Raibert) \\
\texttt{fast} & fast 15 m & SAC fast 0.00: fails & 1/4 / 1/3 / 0/3 & return 6003 (O, Raibert) \\
\texttt{parity} & fast 15 m, return $\ge$ 5020 & SAC fast 0.00: fails its own bound & 0/3 / 2/3 / 1/3 & return 6315 (F, SIMBICON) \\
\texttt{robust} & forward under 4 perturbations & SAC 0.625 at motors $-20\,\%$: fails; SAC-DR clears & 3/3 / 2/3 / 1/3 & return 5563 (O, SIMBICON) \\
\texttt{robust-parity} & brisk, return $\ge$ 5020, 5 siblings & SAC clears by construction; SAC-DR brisk 0.58--0.95: fails & 0/3 / 1/3 / 1/3 & return 6027 (S, SIMBICON) \\
\texttt{robust-noisy} & forward under the suite and noise + delay & SAC 0.05 under noise + delay: fails; SAC-DR (all) misses one sibling at 0.875 & 2/3 / 1/3 / 0/3 & return 5271 (O, SIMBICON) \\
\texttt{push} & forward under 100 N pushes & SAC 0.60: fails; SAC-DR (all) 0.95: clears & 0/3 / 0/3 / 0/3 & open; best 0.675 on push\_v1 (O, SIMBICON) \\
\texttt{robust-noisy-push} & suite, noise + delay and pushes & SAC fails; SAC-DR (all) misses one sibling at 0.875 & 0/3 / 0/3 / 0/3 & open; best 0.825 on noisy\_delay\_v1 (O, SIMBICON) \\
\addlinespace[2pt]\multicolumn{5}{@{}l}{\emph{Lunar lander}} \\
\texttt{standard} & land, achievement 0.9 & PPO 263.9: clears & 3/3 / 3/3 & return 274 (F, PID/PD) \\
\texttt{robust} & terminal under 3 physics scenarios, return $\ge$ PPO $-$ 5 & PPO 263.9 (bound 258.9) & 3/3 / 3/3 & return 281 (O, other) \\
\bottomrule
\end{tabular}

\caption{Control performance of three LLM agents against the RL solutions. Designed controllers reach
the RL return on every parity task, but only in a minority of sessions; robustness at the RL return is reached on the humanoid only by a designed controller, while
only domain-randomised RL clears the 100~N push. Each row is one task: what it asks for, how the RL
solution fares, passes over sessions per model (O Opus 5.5, F Fable 5.1, S Sonnet 5.5), and the best
designed controller, i.e., the pass with the highest return or, on an open task, the attempt closest to
the gate. SAC-DR: SAC trained with domain randomisation. Generated from the session records.}
\label{tab:boundary}
\end{table}

The paper is structured as follows: Section~\ref{sec:benchmark} defines the benchmark, Section~\ref{sec:protocol} the protocol;
Section~\ref{sec:boundary} compares designed controllers with the RL solutions,
Section~\ref{sec:agents} describes the sessions, and Section~\ref{sec:conclusion}
concludes.

\section{The benchmark}
\label{sec:benchmark}

A task asks for one Python file exporting one function\footnote{The function
\texttt{controller(obs)} in the code.} from the plant's output to an action array, starting from a file
that applies zero action. Two rules hold: (1) no neural networks and no gradient-based learning, so that
the artefact stays a white-box program an engineer can read, and (2) a budget of one hour of wall-clock
time, with a submission below the gate scoring zero. The agent
evaluates its controller on visible seeds and has a parameter tuner that searches a declared
block of parameters; grading runs on held-out seeds in a separate container. Each
task is packaged as a Harbor task \citep{harbor2026,terminalbench2026}, so that any agent harness runs
against it. Since it carries a verifiable scalar reward, the proposed ControlHarbor can also serve as an
environment for post-training language models. The RL solutions of Section~\ref{sec:boundary} and the
coding agent act at different levels: one learns the controller from millions of plant steps, the other
writes the controller as code after a few hundred plant episodes. Both levels are graded by the same
gate using the same seeds. The task
format, the sandbox and the task generator are described in Appendix~\ref{app:harbor},
where Figure~\ref{fig:structure} shows how the parts fit together.

\subsection{The grader}
\label{sec:grader}

\paragraph{Problem statement.}
We consider the plant given by a discrete-time system
\begin{equation*}
  x_{t+1} = f_\theta(x_t, u_t, w_t), \qquad y_t = h_\theta(x_t, v_t), \qquad t = 0,\dots,T-1,
\end{equation*}
where $\theta$ selects the plant model from a finite set $\Theta$, which contains the nominal model
$\theta_0$ and the models of the perturbation suite, and a seed $s$ fixes the initial state and the
noise sequences $\{w_t\}$ and $\{v_t\}$.
A controller is a feedback law $\pi : y \mapsto u$ that does not use a neural network. We denote the
set of such laws by $\Pi$.
The requirements of a plant form an ordered set of levels (a ladder), called rungs, $k = 1,\dots,K$,
e.g.\ \{\texttt{alive}, \texttt{forward}, \texttt{brisk}, \texttt{fast}\} in the biped tasks. Rung $k$
is examined through a pass/fail test $\varphi_k$ on each trajectory $\xi^{\pi,\theta}_s$, e.g.\ for \texttt{brisk}: no fall and at least
8~m forward; we write $\xi \models \varphi_k$ when the trajectory passes, and $\varphi_{k+1} \Rightarrow \varphi_k$. A task fixes a target rung $k$ and asks for
\begin{equation}\label{eq:problem}
  \text{find } \pi \in \Pi \ \text{ s.t. } \
  \Pr_s\big[\xi^{\pi,\theta}_s \models \varphi_k\big] \ge 1-\varepsilon \ \ \forall \theta \in \Theta, \;\;
  J_{\theta_0}(\pi) \ge \rho\, J_{\theta_0}(\pi_{\mathrm{RL}}), \;\;
  C(\pi) = 1,
\end{equation}
with $\varepsilon = 0.1$, $\rho = 0.9$, $J_\theta$ the expected return and $C$ the certificate; the last
two constraints apply only when the task declares them. An RL solution instead maximises
the expected return, on the nominal plant or, under domain randomisation, averaged over $\Theta$; the
requirement in \eqref{eq:problem} must instead hold on each plant in $\Theta$.
The grader estimates each probability on $n$ held-out seeds. If all $n$ pass, the one-sided
Clopper--Pearson bound gives, with confidence $1-\delta$, a failure probability of at most
$1-\delta^{1/n}$: $0.072$ for $n = 40$ and $\delta = 0.05$, but $0.45$ for the five visible seeds the
agent evaluates on.

\paragraph{Implementation.}
Each plant's evaluator rolls the controller out on a list of seeds and returns one record of metrics per
seed, and a generic grader composes the gate \eqref{eq:problem} from these records. The platoon rungs
are \texttt{safe}, \texttt{spacing}, \texttt{string\_1\_35}, \texttt{string\_1\_0}; the locomotion tasks use 40 held-out seeds; siblings are e.g.\ torso mass $+20\,\%$, friction
$-30\,\%$, motors $-20\,\%$, gravity $+10\,\%$ or a 75~N push on the humanoid; and the platoon
certificate requires a frequency sweep of a sinusoidal leader whose amplitude ratio along the string
must not exceed one. The grader writes a reward of one or zero and a record
with the failed seeds and the metrics that explain them, which is what the agent reads on the visible
seeds. One variant moves the tuner into the grader, so that the agent is judged on the structure it
wrote rather than on its numbers. Thresholds come from physics or from a frozen baseline, not from
what an agent reached. A parity bound is $90\,\%$ of the mean return of a public RL solution, or, on the
humanoid, of an RL solution we trained, measured on the grader's seeds. The humanoid robust-parity suite
keeps the perturbations this RL solution still clears on 0.9 of the seeds.

The grader is an empirical test on simulated episodes. It does not give a proof: the platoon
certificate is a frequency sweep, not a proof of string stability, and no rung carries a Lyapunov
argument, a reachability analysis or a test on hardware. Any of these can be added later as a further
component of the gate.

\subsection{Plants and tasks}
\label{sec:plants}

The 13 plants are listed with their ladders in Table~\ref{tab:plants} of the
appendix. The difficult ones are a platoon of 20 vehicles with delayed, noisy local sensing, and the
MuJoCo hopper, walker and humanoid. Within a plant, the 44 tasks vary
the dimensions named in the introduction: on the platoon, delay and noise, heterogeneity, a strict
string-stability rung and the certificate; on the lander, 14 physics scenarios; on the bipeds, the
speed tier, the physics suite, pushes, the parity bound and a tool ablation without the tuner.

\section{Protocol}
\label{sec:protocol}

Three models were run as coding agents under one protocol: Claude Sonnet 5.5, Claude Opus 5.5 and
Claude Fable 5.1,\footnote{Claude Fable 5.2 was released on 9 November 2026, after the sessions reported here were run.} each through the Claude Code harness at its highest reasoning effort, without network access, launched by
Harbor. A session is one container, one instruction, one hour from the agent's first command and one
graded file. The instruction states the budget, a score of zero for any submission below the
gate, and the command to read the clock. Three sessions were run per task and model, so $2/3$ means two of them passed. Opus and Fable ran on all 44 tasks and Sonnet on
the 22 MuJoCo tasks, 330 sessions in all, each sandbox limited to two cores and 2~GB. 

The harness records each session's full stream of messages and tool calls, its token usage and timings,
the final file and the grader's record; Table~\ref{tab:sessions} of the appendix gives the medians per
model. Sonnet uses almost the whole hour, 52 minutes at the median, Opus and Fable about half, and 44
sessions, all passes on the light plants, ended within 10 minutes because the agent had a passing
controller and stopped. A session runs the official grader on the visible seeds two or three times and
writes and runs far more Python of its own; the tuner was invoked in a third of Fable's sessions and a
seventh of the others'. A verifier timeout ended 14 sessions, all of them sampling
model-predictive controllers simulating the plant inside the loop.

Two sources of evidence sit outside the protocol. The RL solutions of Section~\ref{sec:boundary}
are single checkpoints, the public SB3-zoo solutions for the hopper, the solutions for the walker and the lander 
trained by us for the humanoid, all measured on the grader's seeds. We have trajectories of the same plants
through a different harness and seven model families, among them GPT, Gemini and Sol, had no clock and
mostly one session per task; it is reported in Appendix~\ref{app:cursor} as context, not as evidence in this paper.

\section{Designed controllers against RL solutions}
\label{sec:boundary}

\subsection{References}
\label{sec:references}

Four plants carry an RL solution: the public SB3-zoo TD3 solution for the hopper and the walker
\citep{raffin2021sb3,fujimoto2018td3}, the zoo PPO solution for the lunar lander \citep{schulman2017ppo},
and for the humanoid, where the public RL solutions do not transfer to the current simulator, a SAC solution
\citep{haarnoja2018sac} trained by us for 1.8~M steps on the nominal plant. All are single checkpoints of
a single seed, measured on the grader's held-out seeds: TD3 returns 3532 on the hopper and 4701 on the
walker, PPO 263.9 on the lander, SAC 5579 on the humanoid. The parity bounds are $90\,\%$ of the
hopper, walker and humanoid returns: $3179$, $4231$ and $5020$. A nominal RL solution is not a fair reference for a robustness gate, so we also trained SAC with
domain randomisation, written SAC-DR below and graded with the same verifier. The hopper and the walker
each have one such solution, with the physics perturbations redrawn every episode. Two SAC-DR solutions were trained for the humanoid: one under physics perturbations and pushes,
and one, written SAC-DR over everything, under those perturbations together with noise and delay.

\subsection{Comparison with the RL solution}
\label{sec:parity}

Table~\ref{tab:boundary} gives, per tier of the four plants, what the gate asks, what the RL solution does
on it, the passes per model and the best controller designed in session. A nominal RL solution clears the forward and brisk tiers. The agents clear these tiers in most
sessions on the hopper and the walker, and in about half of the sessions on the humanoid. On the parity
tiers the gate asks for the RL solution's speed and for $90\,\%$ of its return, and a designed controller
cleared it on every plant. The best returns are 3450 on the hopper, $98\,\%$ of TD3's 3532; 5265 on the
walker, above TD3's 4701; and 6315 on the humanoid, above SAC's 5579, at the fast speed SAC never
reaches. The passes are a minority, one to three of nine sessions per task, from different models on
different plants. The lunar lander shows no boundary: every session cleared both its gates at PPO's return.

What passes differs by plant. On the humanoid every parity pass is a textbook design, SIMBICON, a
Raibert-style foot-placement walker or a capture-point walker, with gains set by hand or by the tuner.
On the hopper and the walker the parity passes are a linear feedback law on the state and a gait clock,
mirror-symmetric across the legs, whose gains the agent found by a black-box search, e.g., CMA-ES or random
search, that it programmed itself. The controller is close to the static linear policies of
\citep{mania2018ars}. Removing the
tuner from the image changes none of this: without it the agents passed the three parity tasks 8 times
in 27 sessions against 7 with it.

\subsection{Robustness}
\label{sec:robustness}

Every parity pass was re-run on 16 public seeds under the perturbations of the robust tiers, and
on the humanoid under the 75~N push. Table~\ref{tab:parityrobust} gives the rate at which the controller
stayed upright and moved forward. A controller designed for speed alone falls under those perturbations.
The searched linear policies of the hopper and the walker collapse under almost every perturbation, and
the humanoid SIMBICON and foot-placement walkers fall under weaker motors and under the push.

The humanoid robust-parity tier asks for speed and robustness together. SIMBICON solved it, at the
return of the SAC solution, with every perturbation above 0.94, the push included. On the walker the
previous oracle, a SIMBICON at $69\,\%$ of TD3's return, survives all four perturbations, while the
parity pass survives none.

SAC-DR shows the same split. On the walker it clears the robust tier that TD3 fails, at a return below
the parity bound. On the hopper it clears nothing, after one million steps and after two million; the
designed controllers clear that tier in six of nine sessions, and the brisk-robust tier in three. On the
humanoid, SAC-DR over physics and pushes clears the robust tier at walking pace and fails the
robust-parity tier that SIMBICON passes. SAC-DR over everything walks at $0.3$--$0.4$~m/s under noise and
delay, and misses the robust-noisy gate on one sibling at 0.875, a gate the agents reached in three of
nine sessions. Designed controllers reached robustness at walking pace on every plant, and robustness at
the RL solution's speed on the humanoid only.

\subsection{Open tasks}
\label{sec:frontier}

Two tasks on the humanoid remain open for every model. Under 100~N pushes the best designed
controller, a SIMBICON walker by Opus, reaches 0.675 against a gate of 0.9, while SAC-DR over
everything clears that gate at 0.95. Under pushes together with physics changes, noise and delay,
neither side clears the gate. The best designed controller reaches 0.825: it survives the pushes on
97\,\% of the seeds and fails once noise and delay are added. SAC-DR over everything misses the same
gate on one sibling. The parity and robust-parity tiers are passed in fewer than half of the
sessions: 1 of 9 on the hopper, 3 of 9 on the walker and 3 of 9 on the humanoid, and 2 of 9 for
humanoid robust-parity. Forward and robust at walking pace are passed in most sessions.
Designed controllers reach robustness at walking pace on every plant. SAC-DR reaches it on the
walker and the humanoid, and clears nothing on the hopper. Robustness at the RL solution's speed
is reached on the humanoid only, by a designed controller, and SAC-DR does not reach it. The 100~N
push is the reverse: SAC-DR clears it and no designed controller does. Pushes together with physics
changes, noise and delay are missed by both.

\begin{table}[t]
\centering\footnotesize
\setlength{\tabcolsep}{3.5pt}
\begin{tabular}{@{}llrcccccc@{}}
\toprule
task & design (model) & return & nominal & torso & friction & motors & gravity & push \\
\midrule
\addlinespace[2pt]\multicolumn{9}{@{}l}{\texttt{hopper-parity}} \\
 & linear policy tuned by search (O) & 3450 & \textbf{1.00} & 0.00 & 0.00 & 0.00 & 0.00 & -- \\
 & linear policy tuned by search (sol) & 3377 & \textbf{1.00} & \textbf{1.00} & 0.00 & 0.00 & 0.00 & -- \\
 & linear policy tuned by search (O) & 3289 & \textbf{1.00} & 0.19 & 0.00 & 0.00 & 0.00 & -- \\
\addlinespace[2pt]\multicolumn{9}{@{}l}{\texttt{hopper-parity-notune}} \\
 & linear policy tuned by search (O) & 3515 & \textbf{0.94} & 0.00 & 0.00 & 0.00 & 0.06 & -- \\
 & linear policy tuned by search (O) & 3328 & \textbf{0.94} & 0.75 & 0.00 & 0.69 & \textbf{1.00} & -- \\
 & linear policy tuned by search (O) & 3277 & \textbf{0.94} & 0.50 & 0.00 & 0.00 & \textbf{0.94} & -- \\
\addlinespace[2pt]\multicolumn{9}{@{}l}{\texttt{walker2d-parity}} \\
 & linear policy tuned by search (O) & 5265 & \textbf{1.00} & 0.00 & 0.00 & 0.00 & 0.00 & -- \\
 & linear policy tuned by search (O) & 4477 & \textbf{1.00} & 0.00 & 0.00 & 0.00 & 0.00 & -- \\
 & SIMBICON (F) & 4429 & \textbf{1.00} & 0.00 & 0.56 & 0.00 & \textbf{1.00} & -- \\
 & SIMBICON, tuned (previous oracle) & 3254 & \textbf{1.00} & \textbf{1.00} & \textbf{1.00} & \textbf{1.00} & \textbf{1.00} & -- \\
\addlinespace[2pt]\multicolumn{9}{@{}l}{\texttt{walker2d-parity-notune}} \\
 & CPG / clocked oscillator (O) & 4679 & \textbf{1.00} & \textbf{1.00} & 0.81 & 0.00 & 0.88 & -- \\
 & CPG / clocked oscillator (S) & 4371 & \textbf{1.00} & \textbf{1.00} & 0.50 & 0.00 & \textbf{0.94} & -- \\
\addlinespace[2pt]\multicolumn{9}{@{}l}{\texttt{humanoid-parity}} \\
 & SIMBICON (F) & 6315 & \textbf{0.94} & 0.06 & 0.00 & 0.00 & 0.00 & 0.00 \\
 & Raibert / foot placement (S) & 6228 & \textbf{1.00} & 0.00 & 0.00 & 0.00 & 0.88 & 0.00 \\
 & capture point / LIP walking (F) & 6125 & \textbf{1.00} & \textbf{1.00} & 0.88 & 0.00 & 0.12 & 0.00 \\
\addlinespace[2pt]\multicolumn{9}{@{}l}{\texttt{humanoid-parity-notune}} \\
 & SIMBICON (O) & 6351 & \textbf{1.00} & 0.81 & 0.00 & 0.00 & 0.00 & 0.00 \\
 & SIMBICON (F) & 6340 & \textbf{0.94} & \textbf{1.00} & 0.62 & 0.00 & 0.69 & 0.25 \\
 & SIMBICON (F) & 6181 & \textbf{1.00} & 0.25 & 0.00 & 0.00 & 0.00 & 0.00 \\
\addlinespace[2pt]\multicolumn{9}{@{}l}{\texttt{humanoid-robust-parity}} \\
 & SIMBICON (S) & 6027 & \textbf{1.00} & \textbf{1.00} & \textbf{1.00} & \textbf{0.94} & \textbf{1.00} & \textbf{1.00} \\
 & SIMBICON (F) & 5624 & \textbf{1.00} & \textbf{1.00} & \textbf{1.00} & \textbf{1.00} & \textbf{1.00} & \textbf{1.00} \\
\bottomrule
\end{tabular}

\caption{Every parity pass re-run on 16 public seeds under the robust-tier perturbations: forward
rate per scenario (bold where at least 0.9), with the walker's previous oracle for contrast. Torso
$+20\,\%$, friction $-30\,\%$, motors $-20\,\%$, gravity $+10\,\%$, push 75~N (humanoid only).}
\label{tab:parityrobust}
\end{table}

\section{The sessions}
\label{sec:agents}

\paragraph{Design families.} The task puts no constraints on control law. The agent may write any function from the
observation to an action, except for a neural network or gradient-based learning. Each submitted file was
read once afterwards by a language model and assigned a design family.

The files are textbook methods for the plant. On the consensus, formation and attitude plants the
controllers are the standard consensus protocol, 18 of 18 passing; on the pendulum and the cart-link an
energy-based swing-up with a linear-quadratic catch, 12 of 12; on the arm and the platoon feedback
linearisation, 12 of 12, and on the platoon the string-stable cooperative cruise controller, 47 of 53;
on the crowd, barrier functions and social forces, 9 of 9. On these plants the agents stay with that one
method.

On the hopper, the walker and the humanoid the method changes from session to session. In 199 sessions
SIMBICON was written 71 times and passed 43, Raibert-style foot placement 42 and 20, a clocked pattern
generator 29 and 8, a linear policy with searched gains 13 and 7, capture-point walking 7 and 4, and
sampling model-predictive control 14 times with no pass, each stopped by the verifier's budget. When a
design stalled, the agent switched to another of these. No submitted controller fell outside them.
Of the 14 model-predictive attempts, 11 are Opus's; most SIMBICON passes are Fable's and
Sonnet's.

\paragraph{Guarantees.} No task asks the agent for a proof and the grader checks none. On the light plants
the chosen families carry the textbook guarantees and the agents argue from them: Laplacian feedback on the ring, a linear-quadratic catch under an energy-shaping swing-up, platoon gains set to the
string-stability condition the certificate measures. That guarantee applies to the idealised model.
The grader accepts or rejects the submitted file by the held-out seeds. On the locomotion plants no such argument exists: stability results cover particular gaits on reduced models, the Raibert hopper
\citep{koditschek1991hopping} and hybrid zero dynamics \citep{westervelt2003hzd}, not SIMBICON,
foot-placement heuristics or searched linear policies or RL solutions, and none was expected.

\paragraph{How sessions fail.} Of the 115 sessions that did not pass, 57 fell, i.e.\ cleared the rung on
fewer than half of the held-out seeds, 19 were marginal, 20 cleared the nominal plant and missed one
sibling scenario, 5 fell short of the return bound only, and 14 timed out. In 42 of these 115 the agent's
last run on the five visible seeds reached 0.9, and the 40 held-out seeds did not. Five seeds are too
few to tell a gait that stands on every seed from one that stands on only some of them. The passes are the short sessions: a passing session ran the grader three times, 23
programs of its own, and 25 minutes at the median; a failing one, twice, 40 programs and 51 minutes. The
agents do not stop when they are losing; they search until the hour ends. The tuner made no difference,
45 passes in the 71 sessions that invoked it against 170 in the 259 that did not.

\paragraph{Carrying a controller forward.} The humanoid brisk and fast tiers were run twice. One
task starts from the zero-torque file. The other starts from a controller that passed the tier below.
From the zero file, brisk passed in 4 of 9 sessions and fast in 2 of 10. From the previous controller,
brisk passed in 7 of 8 and fast in 7 of 9. That file was not written by the model that continued the
task. Every other session in this paper starts from the zero file, and the sessions do not see one
another.

\paragraph{A plant outside the sweep.} We also experimented with a plant not commonly seen in the control 
textbook, a hanging two-segment pneumatic soft arm with creep, hysteresis and pressure lag, seen by a tip
camera, in three tasks: reach and hold, the same under four perturbations, and a circle to trace
(Appendix~\ref{app:softarm}). All 27 sessions passed, with the
one architecture of the soft-robotics literature: the nominal static map inverted for a feedforward and the
camera error fed back through a damped inverse of its Jacobian. The robust siblings made seven of nine
sessions identify a gain per chamber, three circle sessions added repetitive control over the lap, and six of
nine reached the fine rung our oracle had not. The recipe exists in soft robotics, and the agents found it.

\section{Conclusion}
\label{sec:conclusion}

We set out to test what current frontier coding agents can do when asked to design a controller. For that
test we propose a Harbor task format: 44 tasks over 13 plants, so that any agent harness
runs against them unchanged. A task grades one file as control engineers grade, on held-out seeds,
against a ladder of named requirements, under perturbation suites and certificates, within an hour. 
Under one protocol, three frontier models designed controllers that reach
the return of an RL solution on every plant where the gate asks for it, even though in a minority of sessions; the
controllers are programs of tens to hundreds of lines that an engineer can read. On the four plants
that carry an RL solution, one evaluation of the reference controller is 2 to 40 times faster than one
evaluation of the RL network, measured on CPU. Robustness appeared
only when the specification demanded it, and it was then within reach of a textbook design, on the
humanoid at the RL solution's own return, where a domain-randomised RL solution falls short; designs
whose gains came from in-session search were as brittle as the RL solutions they match. Two tasks
remain open for the designed route: the humanoid under 100~N pushes, and the same plant under those
pushes together with torso mass $+20\,\%$, friction $-30\,\%$, motors $-20\,\%$, gravity $+10\,\%$,
measurement noise and one-step command delay. The second remains open for the RL route as well.

The limits of the study are those stated along the way. The three models are from one vendor, and the
earlier cross-vendor sweep is not comparable. Three sessions per cell give counts, not rates. The grader
is an empirical test on simulated episodes, below an industrial certificate (e.g., formal verification), and the RL solutions
are single checkpoints of single seeds. 

Some future directions follow directly. The first is whether the open tasks mark a limit of capability,
or of one hour from a zero controller: the handoff tiers showed that one controller carried forward doubles the pass count, and a
board of prior attempts, built from visible-seed results so that nothing held out flows forward, is the
natural next experiment; if the push tasks stay open under accumulated attempts, the frontier is a hard
one. The second is to replace some of the present components with the methods used in industrial
practice, such as reachability and formal verification.

\subsection*{Generative AI Use}
Generative AI tools were used for coding and experimental implementation, and for writing and editing.
The coding agents evaluated in this paper are the object of study and are described in
Section~\ref{sec:protocol}. The authors reviewed and verified all AI-assisted content and take
responsibility for the final content of the paper.

\bibliography{refs}

\appendix
\section{The Harbor task format}
\label{app:harbor}

\subsection{Two loops: the plant and the task}
\label{sec:two-loops}

The benchmark proposed here tests whether a frontier coding agent can design a controller.
It has two loops, and we keep them apart by name. The inner one is the \emph{plant}, the same plant an RL solution is trained on: each action
returns an observation and a reward, over a fixed horizon. The step size belongs to the plant. On the MuJoCo plants one step is
8 to 15 milliseconds of simulated time. The RL solutions are trained in this loop, and every
controller in this paper is also evaluated in it. The outer one is the \emph{task}, a Harbor task \citep{harbor2026,terminalbench2026}:
a container with the plant installed, an instruction, and a verifier.
Figure~\ref{fig:task-folder} shows the folder of one such task; every task is generated with this
layout.

\begin{figure}[!ht]
\centering
\begin{tikzpicture}[font=\scriptsize, >=latex,
  box/.style={draw, rounded corners=2pt, align=center, inner sep=3pt, minimum height=1.35cm}]
\node[box, text width=2.9cm] (task) at (1.6,0)
  {\textbf{Task}\\ instruction, starter file applying zero action, plant $f_\theta$ with perturbation
   suite $\Theta$, rung ladder $\varphi_1,\dots,\varphi_K$};
\node[box, text width=3.3cm] (agent) at (5.75,0)
  {\textbf{Agent session} (sandbox)\\ model run by a harness, terminal, visible seeds, parameter tuner,
   limited time budget};
\node[box, text width=1.9cm] (ctrl) at (9.15,0)
  {\textbf{Controller}\\ one Python file, \texttt{controller}\allowbreak\texttt{(obs)}};
\node[box, text width=3.5cm] (grader) at (12.85,0)
  {\textbf{Grader} (separate container)\\ held-out seeds on $\theta_0$ and every $\theta \in \Theta$,
   Clopper--Pearson bound, certificate, return bound from RL};
\node[box, text width=2.7cm, minimum height=0.8cm] (rl) at (9.15,-1.75)
  {\textbf{RL solutions}\\ SAC, TD3, PPO; SAC-DR};
\node[box, text width=3.0cm, minimum height=0.8cm] (out) at (13.35,-1.75)
  {\textbf{Grade}\\ pass/fail per rung $=$ scalar reward};
\draw[->] (task) -- (agent);
\draw[->] (agent) -- (ctrl);
\draw[->] (ctrl) -- (grader);
\draw[->] (rl.east) -| ([xshift=-1.25cm]grader.south);
\draw[->] ([xshift=0.5cm]grader.south) -- ([xshift=0.5cm]grader.south |- out.north);
\end{tikzpicture}
\caption{ControlHarbor separates writing the controller from grading it: the agent sees the task and
visible seeds only, and a separate grader runs the submitted file on held-out seeds on the nominal and
every perturbed model. RL solutions are graded by the same grader on the same seeds, and the grade
doubles as a reward for post-training.}
\label{fig:structure}
\end{figure}
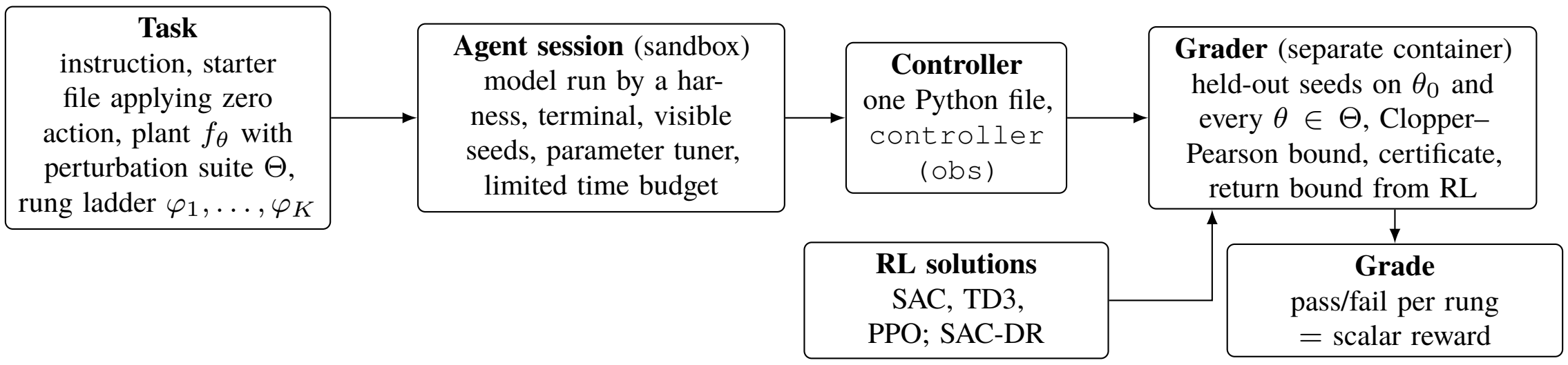

\begin{figure}[!ht]
\centering
\begin{tikzpicture}[y=-0.44cm, x=1cm,
    fn/.style={font=\ttfamily\small, anchor=west, inner sep=1.5pt},
    role/.style={font=\small\itshape, anchor=west, inner sep=1.5pt}]
  \node[fn] at (0,0) {hopper-forward/};
  \node[fn] (md) at (0.55,1) {instruction.md};
  \node[role] at (6.15,1) {what the agent reads};
  \node[fn] (toml) at (0.55,2) {task.toml};
  \node[role] at (6.15,2) {budget, artefact, separate verifier};
  \node[fn] (env) at (0.55,3) {environment/};
  \node[role] at (6.15,3) {the agent's container};
  \node[fn] (docker) at (1.15,4) {Dockerfile};
  \node[role] at (6.15,4) {plant installed, references stripped};
  \node[fn] (starter) at (1.15,5) {starter\_controller.py};
  \node[role] at (6.15,5) {zero action, copied to the work file};
  \node[fn] (sol) at (0.55,6) {solution/};
  \node[role] at (6.15,6) {the oracle};
  \node[fn] (solve) at (1.15,7) {solve.sh};
  \node[role] at (6.15,7) {installs the reference controller};
  \node[fn] (ref) at (1.15,8) {controller.py};
  \node[role] at (6.15,8) {reference controller};
  \node[fn] (tests) at (0.55,9) {tests/};
  \node[role] at (6.15,9) {the verifier};
  \node[fn] (testsh) at (1.15,10) {test.sh};
  \node[role] at (6.15,10) {held-out seeds, fresh container};
  \draw[black!55] (0.22,0.45) -- (0.22,9) -- (md.west |- 0,9);
  \foreach \yy in {1,2,3,6} \draw[black!55] (0.22,\yy) -- (0.5,\yy);
  \draw[black!55] (0.82,3.45) -- (0.82,5) -- (docker.west |- 0,5);
  \draw[black!55] (0.82,4) -- (1.1,4);
  \draw[black!55] (0.82,6.45) -- (0.82,8) -- (solve.west |- 0,8);
  \draw[black!55] (0.82,7) -- (1.1,7);
  \draw[black!55] (0.82,9.45) -- (0.82,10) -- (testsh.west);
\end{tikzpicture}
\caption{The folder of one Harbor task, the hopper forward gate. Every task is generated with
this layout. The agent reads the instruction and edits the starter inside the container. The
verifier grades held-out seeds in a fresh container, and the solution directory is the oracle.}
\label{fig:task-folder}
\end{figure}

Here the agent is a language model
operating a terminal, its actions are shell commands and file edits, an episode is a session with a
wall-clock budget, and the reward is computed once, at the end, on seeds the agent never sees. Nothing
is learned inside this loop: the session trains no neural network and runs no reinforcement-learning
algorithm. It leaves a controller file behind. The RL solution of Section~\ref{sec:boundary} and the coding
agent thus act at different levels, one learning the controller from millions of plant steps, the other
writing it as code after a few hundred plant episodes, and both are graded by the same inner-loop gate
on the same seeds. A Harbor task carries a verifiable scalar reward, thus
 the tasks can also be used as reinforcement learning environments for post-training of language models.

\subsection{The artefact and the sandbox}
\label{sec:contract}

The artefact is one Python file exporting one function, \texttt{controller(obs)}, from the plant's
output to an action array, and the starting file applies zero action. The instruction states the plant, 
its observation and action, the ladder of requirements and the
grader's flags, and imposes two rules: no neural networks and no gradient-based learning, so that the
artefact stays a program an engineer can read, and a budget of one hour from the agent's first command,
with a submission below the gate scoring zero. The sandbox is a container with the plant package
installed and no network; the hidden reference controllers are stripped from the image. The agent has
the evaluator itself as a tool, on visible seeds disjoint from the held-out ones, and a parameter tuner
that searches a declared \texttt{PARAMS} block with a declared budget. It has no clock other than the
shell and is told to run \texttt{date}, which gives it a precise reading of the time. Grading runs in a separate container on the submitted file alone, with a
budget of 15 minutes.

\subsection{Registry and generated tasks}
\label{sec:registry}

Everything specific to a plant, its configuration, evaluator, rungs, failure keys, certificates
and container packages, is declared in one registry; the grader, the tuner and the task generator know
no plant by name. Tasks are generated from one specification, so that the flags the agent reads and the
flags the grader runs cannot drift apart, and a lint checks that no held-out seed appears outside the
grader script. Each task ships an oracle, the best known controller, and passes two fixtures before an
agent sees it: the oracle scores one, the zero-action starter zero; a task whose oracle fails its own
gate is declared open.

\section{Every task under the protocol}
\label{app:tasks}
Table~\ref{tab:tasks} lists every task with its passes per model under the protocol of
Section~\ref{sec:protocol} and the best controller designed in session; the gates are those of
Section~\ref{sec:benchmark} and the task metadata in the repository.

\begin{table}[p]
\centering\scriptsize
\resizebox{\linewidth}{!}{
\setlength{\tabcolsep}{3pt}
\renewcommand{\arraystretch}{0.92}
\begin{tabular}{@{}lc>{\raggedright\arraybackslash}p{4.3cm}>{\raggedright\arraybackslash}p{3.0cm}@{}}
\toprule
task & Opus / Fable / Sonnet & best designed in session & family of the best \\
\midrule
\texttt{humanoid-brisk} & 2/3 / 1/3 / 1/3 & return 5761 (O) & Raibert \\
\texttt{humanoid-brisk-handoff} & 3/3 / 1/2 / 3/3 & return 5929 (S) & SIMBICON \\
\texttt{humanoid-fast} & 1/4 / 1/3 / 0/3 & return 6003 (O) & Raibert \\
\texttt{humanoid-fast-handoff} & 2/3 / 3/3 / 2/3 & return 6080 (F) & SIMBICON \\
\texttt{humanoid-forward} & 2/3 / 2/3 / 1/3 & return 5522 (F) & Raibert \\
\texttt{humanoid-parity} & 0/3 / 2/3 / 1/3 & return 6315 (F) & SIMBICON \\
\texttt{humanoid-parity-notune} & 1/3 / 2/3 / 0/3 & return 6351 (O) & SIMBICON \\
\texttt{humanoid-push} & 0/3 / 0/3 / 0/3 & open; best 0.675 on push\_v1 (O) & SIMBICON \\
\texttt{humanoid-robust} & 3/3 / 2/3 / 1/3 & return 5563 (O) & SIMBICON \\
\texttt{humanoid-robust-noisy} & 2/3 / 1/3 / 0/3 & return 5271 (O) & SIMBICON \\
\texttt{humanoid-robust-noisy-push} & 0/3 / 0/3 / 0/3 & open; best 0.825 on noisy\_delay\_v1 (O) & SIMBICON \\
\texttt{humanoid-robust-parity} & 0/3 / 1/3 / 1/3 & return 6027 (S) & SIMBICON \\
\addlinespace[3pt]
\texttt{hopper-brisk-robust} & 1/3 / 0/3 / 2/3 & return 2652 (O) & other \\
\texttt{hopper-forward} & 3/3 / 3/3 / 2/4 & return 2605 (O) & Raibert \\
\texttt{hopper-forward-v} & 3/3 / 2/3 / 1/3 & return 2203 (F) & other \\
\texttt{hopper-parity} & 1/3 / 0/3 / 0/3 & return 3450 (O) & linear policy (search) \\
\texttt{hopper-parity-notune} & 3/3 / 0/3 / 0/3 & return 3515 (O) & linear policy (search) \\
\texttt{hopper-robust} & 3/3 / 2/3 / 1/3 & return 2145 (F) & other \\
\addlinespace[3pt]
\texttt{walker2d-forward} & 3/3 / 3/3 / 3/3 & return 2494 (S) & CPG \\
\texttt{walker2d-parity} & 2/3 / 1/3 / 0/3 & return 5265 (O) & linear policy (search) \\
\texttt{walker2d-parity-notune} & 1/3 / 0/3 / 1/3 & return 4679 (O) & CPG \\
\texttt{walker2d-robust} & 3/3 / 3/3 / 3/3 & return 2588 (F) & SIMBICON \\
\addlinespace[3pt]
\texttt{platoon-combo-20} & 3/3 / 2/3 & return -167 (F) & CACC \\
\texttt{platoon-combo-20-noisy} & 1/3 / 1/3 & return -515 (O) & CACC \\
\texttt{platoon-combo-20-strict} & 3/3 / 2/3 & return -119 (O) & CACC \\
\texttt{platoon-delay-noise-20} & 3/3 / 3/3 & return -62 (F) & CACC \\
\texttt{platoon-delay-noise-20-certified} & 2/3 / 2/3 & return -64 (O) & CACC \\
\texttt{platoon-delay-noise-20-peak} & 3/3 / 3/3 & return -62 (F) & CACC \\
\texttt{platoon-delay-noise-20-staged} & 3/3 / 3/3 & return -50 (O) & CACC \\
\texttt{platoon-delay-noise-20-strict} & 3/3 / 3/3 & return -62 (F) & CACC \\
\texttt{platoon-delay-noise-20-tuned} & 3/3 / 3/3 & return -79 (F) & CACC \\
\texttt{platoon-delay-noise-20-v} & 3/3 / 3/3 & return -47 (O) & CACC \\
\addlinespace[3pt]
\texttt{crowd-bottleneck-56} & 3/3 / 3/3 & return 442 (O) & other \\
\texttt{crowd-random-exit-32} & 3/3 / 3/3 & return 861 (F) & sampling MPC \\
\addlinespace[3pt]
\texttt{lunar-lander-robust} & 3/3 / 3/3 & return 281 (O) & other \\
\texttt{lunar-lander-standard} & 3/3 / 3/3 & return 274 (F) & PID/PD \\
\addlinespace[3pt]
\texttt{arm-pickup-lag} & 3/3 / 3/3 & return -5 (F) & feedback linearization \\
\texttt{arm-pickup-lag-v} & 2/2 / 3/3 & return -5 (O) & feedback linearization \\
\addlinespace[3pt]
\texttt{cart-link-swingup} & 3/3 / 3/3 & return 564 (O) & energy-based swing-up + LQR catch \\
\addlinespace[3pt]
\texttt{pendulum-swingup} & 3/3 / 3/3 & return 14 (F) & energy-based swing-up + LQR catch \\
\addlinespace[3pt]
\texttt{consensus-ring-8} & 3/3 / 3/3 & return -10 (O) & consensus protocol \\
\addlinespace[3pt]
\texttt{second-order-consensus-ring-8} & 3/3 / 3/3 & return -72 (F) & PID/PD \\
\addlinespace[3pt]
\texttt{formation-circle-8} & 3/3 / 3/3 & return -628 (F) & PID/PD \\
\addlinespace[3pt]
\texttt{so3-sync-ring-6} & 3/3 / 3/3 & return -1 (F) & consensus protocol \\
\bottomrule
\end{tabular}
}
\caption{All tasks under the protocol: passes over sessions per model, the best designed controller
(return of the best pass, or the closest attempt on an open task, with the model) and the design family
of that controller as read from its code. Generated from the session records.}
\label{tab:tasks}
\end{table}

\section{Plants and sessions}
\label{app:plants}
Table~\ref{tab:plants} lists the 13 plants with their rung ladders and certificates,
and Table~\ref{tab:sessions} the sessions per model under the protocol.

\begin{table}[h]
\centering\footnotesize
\setlength{\tabcolsep}{3pt}
\begin{tabular}{@{}>{\raggedright\arraybackslash}p{2.7cm}>{\raggedright\arraybackslash}p{3.4cm}>{\raggedright\arraybackslash}p{4.4cm}>{\raggedright\arraybackslash}p{2.4cm}c@{}}
\toprule
plant & what & rung ladder (gate rungs in order) & certificates & tasks \\
\midrule
Consensus & single integrators on a ring graph & \texttt{converged} & -- & 1 \\
Second-order consensus & double integrators on a ring graph & \texttt{converged} & -- & 1 \\
Formation control & planar double integrators, collision term & \texttt{converged} & -- & 1 \\
SO(3) synchronisation & rigid bodies on SO(3), ring graph & \texttt{converged} & -- & 1 \\
Inverted pendulum & torque-limited swing-up & \texttt{achievement} $\to$ \texttt{terminal} & -- & 1 \\
Cart-link swing-up & cart and link, bounded track & \texttt{achievement} $\to$ \texttt{terminal} & -- & 1 \\
Crowd evacuation & decentralised pedestrians, narrow exit & \texttt{achievement} $\to$ \texttt{terminal} & -- & 2 \\
Vehicle platoon & 20 followers, delayed noisy local sensing & \texttt{safe} $\to$ \texttt{spacing} $\to$ \texttt{string\_1\_35} $\to$ \texttt{string\_1\_0} & string-frequency & 10 \\
Lunar lander & Box2D lander, 14 physics scenarios & \texttt{achievement} $\to$ \texttt{terminal} & -- & 2 \\
Two-link arm & planar arm, lag, unknown payload & \texttt{safe} $\to$ \texttt{reaches} $\to$ \texttt{settles} $\to$ \texttt{tracks} & -- & 2 \\
Hopper & MuJoCo Hopper-v5 & \texttt{alive} $\to$ \texttt{forward} $\to$ \texttt{brisk} $\to$ \texttt{fast} & -- & 6 \\
Walker2d & MuJoCo Walker2d-v5 & \texttt{alive} $\to$ \texttt{forward} $\to$ \texttt{brisk} $\to$ \texttt{fast} & -- & 4 \\
Humanoid & MuJoCo Humanoid-v5 & \texttt{alive} $\to$ \texttt{forward} $\to$ \texttt{brisk} $\to$ \texttt{fast} & -- & 12 \\
\midrule
all & 13 plants & & & 44 \\
\bottomrule
\end{tabular}

\caption{The 13 plants, their rung ladders and certificates, and the number of
tasks on each. Generated from the registry.}
\label{tab:plants}
\end{table}

\begin{table}[h]
\centering\footnotesize
\setlength{\tabcolsep}{4pt}
\begin{tabular}{lrrrrrrr}
\toprule
model & sessions & tasks & passes & minutes & prompt tokens & tool calls & grader runs \\
\midrule
Opus 5.5 & 132 & 44 & 101 & 29 & 2.5\,M & 41 & 2 \\
Fable 5.1 & 131 & 44 & 91 & 31 & 3.5\,M & 40 & 3 \\
Sonnet 5.5 & 67 & 22 & 23 & 52 & 3.2\,M & 55 & 3 \\
\bottomrule
\end{tabular}

\caption{Sessions per model under the protocol. Minutes, prompt tokens (cached tokens included), tool
calls and runs of the official grader on the visible seeds are medians over the model's sessions.}
\label{tab:sessions}
\end{table}

\section{Protocol details}
\label{app:protocol}
Every instruction carries the same closing rules, reproduced here from the hopper parity task with the
seeds omitted: ``The grader uses different seeds than the ones above and the same thresholds and flags.
There is no network access. Your budget is one session of at most 60 minutes of wall-clock time,
counted from your first command; a submission below the gate scores zero, so keep a passing version in
\texttt{/app/work/controller.py} at all times once you have one. You have no other clock than the
shell: run \texttt{date} with your first command and again whenever you decide whether to continue;
the number of turns you have made says nothing about the time used.'' The instruction opens with the
two design rules, no training of neural networks and no gradient-based learning, and lists the plant,
the observation and action, the ladder with its thresholds, the perturbation scenarios, the exact
grader command on the visible seeds, and the tuner's interface.

Harbor launched each session with the Claude Code agent at the highest reasoning effort, three
sessions per task and model, five to six concurrently, with the agent's own setup allowed five times
Harbor's default; the sandbox is the task's image with two cores, 2~GB of memory, no network and a
one-hour agent timeout, and the verifier runs in a separate container from the same image with a
15-minute timeout on the submitted file alone. The verifier flags of every task are in its grader
script in the released repository; the held-out seed blocks live only there. Tokens and minutes per
session are those of Table~\ref{tab:sessions}; the longest session used 24~M prompt tokens, most of
them cached.

\section{Earlier sweep under a different harness}
\label{app:cursor}
Before the protocol of Section~\ref{sec:protocol} was fixed, the same tasks were run through the Cursor CLI
harness (2026-10-21 to 23) with seven model families. These sessions are not comparable with the main
results: a different harness, mostly one session per task, and no clock in the instruction, so a model that
believed its hour was over stopped early (Fable in particular). They are reported here because they are the
only runs with models from other vendors. Nine platoon trials graded under an earlier version of their task
are left out. Table~\ref{tab:cursor} gives passes over sessions on the grader's held-out seeds; a dash means
the model was not run on that task. Every family solved the light plants; the locomotion plants
separated them, with the humanoid at 1 of 3, 0 of 4 and 0 of 1 and the hopper's ladder passed by Sol
and failed by Gemini; GPT and Composer were run on the light plants and a few harder ones only, where
GPT failed the lander and the platoon. None of the parity, robust or ablation tasks existed at the time.

\begin{table}[h]
\centering\footnotesize
\setlength{\tabcolsep}{3.5pt}
\begin{tabular}{lccccccc}
\toprule
task & Opus 5.5 & Fable 5.1 & Sonnet 5.5 & Sol & Gemini & GPT & Composer \\
\midrule
\addlinespace[2pt]\multicolumn{8}{l}{\emph{Humanoid}} \\
\texttt{humanoid-forward} & 1/3$^{\dagger}$ & 0/4$^{\dagger}$ & -- & -- & 0/1 & -- & -- \\
\addlinespace[2pt]\multicolumn{8}{l}{\emph{Hopper}} \\
\texttt{hopper-brisk-robust} & 0/1 & 0/1 & -- & -- & -- & -- & -- \\
\texttt{hopper-forward} & 3/3 & 2/3 & -- & 3/3 & 0/3 & -- & -- \\
\texttt{hopper-forward-v} & 2/2 & 0/3 & -- & 0/1 & -- & -- & -- \\
\texttt{hopper-parity} & 1/2 & 0/2 & -- & 1/1 & -- & -- & -- \\
\texttt{hopper-robust} & 1/2 & 0/2 & -- & 1/2 & -- & -- & -- \\
\addlinespace[2pt]\multicolumn{8}{l}{\emph{Walker2d}} \\
\texttt{walker2d-forward} & 4/4$^{\dagger}$ & 2/3$^{\dagger}$ & -- & 2/2 & 2/2 & -- & -- \\
\addlinespace[2pt]\multicolumn{8}{l}{\emph{Vehicle platoon}} \\
\texttt{platoon-combo-20} & 1/1 & 1/1 & -- & 1/1 & 0/1 & -- & -- \\
\texttt{platoon-delay-noise-20} & 2/2 & 3/3 & 1/1 & 3/3 & 3/3 & 0/1 & -- \\
\texttt{platoon-delay-noise-20-peak} & 2/2 & 2/2 & -- & 1/2 & 2/2 & -- & -- \\
\texttt{platoon-delay-noise-20-staged} & 1/1 & 1/1 & -- & 1/1 & 1/1 & 0/1 & 0/1 \\
\texttt{platoon-delay-noise-20-tuned} & 2/2 & 2/2 & -- & 2/2 & 1/2 & -- & -- \\
\texttt{platoon-delay-noise-20-v} & 2/2 & 2/2 & -- & 2/2 & 0/1 & -- & -- \\
\addlinespace[2pt]\multicolumn{8}{l}{\emph{Crowd evacuation}} \\
\texttt{crowd-bottleneck-56} & 1/1 & 1/1 & -- & 1/1 & 1/1 & -- & -- \\
\texttt{crowd-random-exit-32} & 1/1 & 1/1 & -- & 1/1 & 1/1 & -- & -- \\
\addlinespace[2pt]\multicolumn{8}{l}{\emph{Lunar lander}} \\
\texttt{lunar-lander-robust} & 2/2 & 2/2 & -- & 2/2 & 2/2 & 0/2 & -- \\
\texttt{lunar-lander-standard} & 1/1 & 1/1 & -- & 1/1 & 1/1 & 0/1 & -- \\
\addlinespace[2pt]\multicolumn{8}{l}{\emph{Two-link arm}} \\
\texttt{arm-pickup-lag} & 2/2 & 2/2 & -- & 2/2 & 2/2 & -- & -- \\
\texttt{arm-pickup-lag-v} & 2/2 & 2/2 & -- & 2/2 & 1/2 & -- & -- \\
\addlinespace[2pt]\multicolumn{8}{l}{\emph{Cart-link}} \\
\texttt{cart-link-swingup} & 1/1 & 1/1 & -- & 1/1 & 1/1 & -- & -- \\
\addlinespace[2pt]\multicolumn{8}{l}{\emph{Inverted pendulum}} \\
\texttt{pendulum-swingup} & 1/1 & 1/1 & 1/1 & 1/1 & 1/1 & 1/1 & 1/1 \\
\addlinespace[2pt]\multicolumn{8}{l}{\emph{Consensus}} \\
\texttt{consensus-ring-8} & 1/1 & 1/1 & 1/1 & 1/1 & 1/1 & 1/1 & 1/1 \\
\addlinespace[2pt]\multicolumn{8}{l}{\emph{Second-order consensus}} \\
\texttt{second-order-consensus-ring-8} & 1/1 & 1/1 & 1/1 & 1/1 & 1/1 & 1/1 & 1/1 \\
\addlinespace[2pt]\multicolumn{8}{l}{\emph{Formation control}} \\
\texttt{formation-circle-8} & 1/1 & 1/1 & 1/1 & 1/1 & 1/1 & 1/1 & 1/1 \\
\addlinespace[2pt]\multicolumn{8}{l}{\emph{SO(3) synchronisation}} \\
\texttt{so3-sync-ring-6} & 1/1 & 1/1 & 2/2 & 1/1 & 1/1 & 1/1 & 1/1 \\
\midrule
all tasks & 37/42 & 30/44 & 7/7 & 32/35 & 23/31 & 5/10 & 5/6 \\
\bottomrule
\end{tabular}

\caption{Cursor-era sweep: passes/sessions per task and model family on the held-out seeds. $\dagger$: at
least one of the sessions had the clock in the instruction; all others ran without it.}
\label{tab:cursor}
\end{table}

\section{A plant outside the sweep: the soft arm}
\label{app:softarm}
The arm is a pseudo-rigid-body chain of 12 links on hinge joints in MuJoCo, 30~cm long and hanging under
gravity, with two segments each driven by an antagonistic pair of pressure chambers, i.e.\ four commands in
$[0,1]$ at 100~kPa full scale. A material layer adds creep (a standard linear solid), co-contraction
stiffening, a Bouc--Wen hysteresis loop on the actuation, a first-order pressure lag of 0.15~s with a slew
limit, and a draw of stiffness, damping, chamber gains and loop width per seed; a tip camera at 25~Hz with
1~mm noise is the only measurement of the shape, and the instruction provides the nominal arm's static tip
map. The tasks and their gates: \texttt{soft-arm-reach}, six targets each held within 10~mm for 2~s (fine);
\texttt{soft-arm-reach-robust}, the same on the nominal plant and under a 5~g payload, a 10~g hook attached
mid-episode, a stiffness drift and one chamber at 75\,\%; \texttt{soft-arm-trace}, a 3~cm circle with a
period of 8~s within 20~mm RMS (coarse), with the fine (10~mm) and fast (5~mm) rungs left open because the
oracle, the map inverse with a PI-D correction on the camera, tracked at 15~mm. The sessions ran under the
protocol of Section~\ref{sec:protocol} after the main sweep and are not counted in its tables.
Table~\ref{tab:softarm} gives the passes, the median session length and the number of sessions that reached
each open rung on the held-out seeds. Every controller inverts the nominal map for a feedforward and corrects
on each camera frame through a damped inverse of the map's Jacobian; no session built a simulator of the arm.
On the robust task seven of nine sessions identify a gain per chamber from the commands that earlier targets
converged to, and the chamber at 75\,\% is the binding sibling at 0.95--1.0. On the circle seven of nine lead
the target to cover the pressure lag, three add a repetitive controller indexed by the phase of the lap, and
one inverts the hysteresis through a model of the loop; the two sessions without a lead are the two lowest
fine rates. The plant, evaluator, tasks and records are in the repository with the rest, and the build record
is \texttt{docs/soft\_arm\_spec.md}.

\begin{table}[h]
\centering\footnotesize
\setlength{\tabcolsep}{4pt}
\begin{tabular}{@{}llccccc@{}}
\toprule
task & gate & Opus / Fable / Sonnet & min (median) & fine $\geq 0.9$ & fast $\geq 0.9$ \\
\midrule
\texttt{soft-arm-reach} & fine & 3/3 / 3/3 / 3/3 & 19 & gate & 1/9 \\
\texttt{soft-arm-reach-robust} & fine & 3/3 / 3/3 / 3/3 & 38 & gate & 1/9 \\
\texttt{soft-arm-trace} & coarse & 3/3 / 3/3 / 3/3 & 29 & 6/9 & 1/9 \\
\bottomrule
\end{tabular}

\caption{The soft arm tasks: gated rung, passes over sessions per model, median session length in minutes,
and sessions of nine that reached each rung above the gate on 0.9 of the held-out seeds. Generated from the
session records.}
\label{tab:softarm}
\end{table}

\FloatBarrier
\section{Reproducibility}
\label{app:repro}
The repository at \url{https://github.com/jq-wei/hl-control} holds the plants and evaluators (\path{src/hl_control/envs},
\texttt{eval}), the registry, grader, tuner and task generator (\path{src/hl_control/harbor}), the
44 tasks and the three soft arm tasks (\texttt{harbor/tasks}, each with instruction, image, grader
script and oracle), the
RL solutions' grades, the designed controllers of every pass with their grades, and the tooling:
the ledger, the per-task boards, the video portal (online at \url{https://jq-wei.github.io/hl-control/}) and the generator of every table in this paper. A
task is generated from one specification by \texttt{hl\_control.harbor.scaffold}, which renders the
instruction and the grader script together and validates the oracle and the zero-action starter; a sweep
is one queue file, one line per task, model list and session count, run by \texttt{scripts/run\_harbor\_queue.sh}
through Harbor; the ledger, the boards, the portal and the tables are regenerated by one script each
from the session records. The session transcripts themselves are retained by the authors.

\end{document}